\documentclass{amsart}
\usepackage{amsmath}
\usepackage{graphicx}
\usepackage{latexsym}
\usepackage{epsfig}

\newtheorem{thm}{Theorem}[section]
\newtheorem{lem}[thm]{Lemma}
\newtheorem{cor}[thm]{Corollary}

\theoremstyle{definition}
\newtheorem{defn}[thm]{Definition}
\newtheorem{remark}[thm]{Remark}
\newtheorem{example}[thm]{Example}

\newcommand{\R}{\mathbb{R}}
\newcommand{\Z}{\mathbb{Z}}
\newcommand{\N}{\mathbb{N}}

\begin{document}

\title[Tensor products of quandles and 1-handles]
{Tensor products of quandles and 1-handles attached to surface-links}
\vspace{0.2in}

\author[Seiichi Kamada]{Seiichi Kamada}
\address{Department of Mathematics,
Osaka University,
Toyonaka, Osaka 560-0043, Japan}
\email{kamada@math.sci.osaka-u.ac.jp}




\begin{abstract}
A quandle is an algebra with two binary operations 
satisfying three conditions which are related to  Reidemeister moves in knot theory.   
In this paper we introduce the notion of the (canonical) 
tensor product of a quandle.  
The tensor product of the knot quandle or the knot symmetric quandle 
of a surface-link in $4$-space can be used to  
classify or construct invariants of $1$-handles attaching to the surface-link.   
We also compute the tensor products for dihedral quandles and their symmetric doubles. 
\end{abstract}

\maketitle


\section{Introduction}

A quandle $X =(X, \ast, \overline{\ast})$ is a set $X$ with two binary operations $\ast$ and $\overline{\ast}$ 
satisfying conditions related to Reidemeister moves in knot theory (cf. \cite{Joyce82, Matveev82} and Section~\ref{sect:quandles}).  
A symmetric quandle is a quandle equipped with an involution called a good involution (\cite{Kamada07, Kamada14A, KO10}).  
For a given quandle $X$, there is a symmetric quandle $(D(X), \rho)$, called the symmetric double of $X$ (Example~\ref{doublecover}). 
Quandles and symmetric quandles are used in knot theory for construction of invariants of links and surface-links. 

In this paper, we define the {\it tensor product} $X \otimes X$ of a quandle $X$.   
We give a method of computing the tensor product $D(X) \otimes D(X)$ of the symmetric double $D(X)$ 
from the tensor product $X \otimes X$ of $X$  (Theorem~\ref{thm:tensordouble}).  
Let $\tau: X \otimes X \to X \otimes X$ be the involution switching the components, i.e., 
$\tau(x, y) = (y,x)$.  When $X$ is equipped with 
a good involution $\rho: X \to X$, we have another involution $\rho$ ($=\rho \otimes \rho$)  $: 
X \otimes X \to X \otimes X$ with $\rho (x,y) = (\rho(x), \rho(y))$. 
We consider the quotient set $X \otimes X / \langle \tau \rangle$ ($X \otimes X / \langle \rho \rangle$
or $X \otimes X/ \langle \tau, \rho \rangle$) of $X \otimes X$ by the equivalence relation generated by $\tau$ (by $\rho$ or by $\tau$ and $\rho$).  
It will be seen that the tensor product $X \otimes X$ itself and 
the quotient set $X \otimes X / \langle \tau \rangle$ (or $X \otimes X/ \langle \tau, \rho \rangle$) of $X \otimes X$ are useful for studying $1$-handles attached to surface-links.  

We compute the tensor product $R_n \otimes R_n$ for the dihedral quandle $R_n$ of every order $n$ 
(Theorems~\ref{thm:dihedral:odd} and \ref{thm:dihedral:even}) in Section~\ref{sect:dihedral}. 
Once we know the elements of $R_n \otimes R_n$, we can obtain the quotient set 
$R_n \otimes R_n/ \langle \tau \rangle$ (Corollaries~\ref{cor:dihedral:odd} and \ref{cor:dihedral:even}). 
For the symmetric double $(D(R_n), \rho)$ of $R_n$, using Theorem~\ref{thm:tensordouble}, 
we can compute  $D(R_n) \otimes D(R_n)$ from $R_n \otimes R_n$ 
(Theorems~\ref{thm:d:dihedral:odd} and \ref{thm:d:dihedral:even}), 
and then we obtain the quotient sets $D(R_n) \otimes D(R_n)/ \langle \tau \rangle$, 
$D(R_n) \otimes D(R_n)/ \langle \rho \rangle$ and $D(R_n) \otimes D(R_n)/
 \langle \tau, \rho \rangle$ in Section~\ref{sect:d:dihedral}.

By a {\it surface-link}, we  mean a closed surface embedded in $\R^4$.  
The surface may be disconnected or non-orientable. When it is orientable and oriented, it is called an {\it oriented} surface-link.   A surface-link is called a {\it surface-knot} when it is connected. 
When $F$ is oriented, the {\it knot quandle} $Q(F)$ is defined (cf. \cite{FennRourke92, Joyce82, Kamada07, Kamada14A, KO10, Matveev82} and Example~\ref{example:knotqdle}).   
For any surface-link $F$, the {\it knot symmetric quandle} $(\widetilde{Q}(F), \rho)$ is defined (cf. \cite{Kamada07, Kamada14A, KO10} and 
Example~\ref{example:knotsymqdle}).  

Using tensor products of quandles, we discuss classification and invariants of $1$-handles attached to surface-links.  

F. Hosakawa and A. Kawauchi \cite{HK} studied unknotted surface-links and surgery along $1$-handles.   They proved that any oriented surface-link can be transformed to an unknotted oriented  surface-link by surgery along oriented $1$-handles.  For non-orientable surface-links, an analoguos result was shown by the author in \cite{Kam89}, that is, any non-orientable surface-link can be transformed to an  unknotted non-orientable surface-link by surgery along $1$-handles.   In this sense, surgery along a $1$-handle is often called an unknotting operation for surface-links.  Surgery along a $1$-handle is also used to construct examples of surface-knots and surface-links.  

J. Boyle \cite{Boyle} classified oriented $1$-handles attached to 
an oriented surface-knot $F$ up to strong equivalence in terms of the knot group $G(F) :=  \pi_1(\R^4 \setminus F)$ and its peripheral subgroup.  They are in one-to-one correspondence to double cosets of the knot group $G(F)$ by the peripheral subgroup.    
The author \cite{Kamada14B} extended Boyle's argument to the case where $F$ is a surface-knot which may be unoriented or non-orientable.  However the arguments in \cite{Boyle} and \cite{Kamada14B} using knot groups are not applied directly to surface-links which are not connected.  
Even for  surface-knots,  classifications of $1$-handles stated in terms of the knot groups given in \cite{Kamada14B} are not easy to discuss.  
The argument using tensor products discussed in this paper is simple and applicable to any surface-link which may be disconnected or non-orientable.  

For a surface-link $F$, let 
$\mathcal{H}^{\rm s}(F)$ (or $\mathcal{H}^{\rm w}(F)$)  denote the 
set of strong (or weak) equivalence classes of $1$-handles attached to $F$.   
(See Section~\ref{sect:handle} for the definition of strong/weak equivalence.)

\begin{thm}\label{thm:mainA}
Let $F$ be a surface-link and $(\widetilde{Q}(F), \rho)$ the knot symmetric quandle of $F$.  
There are bijections
\begin{equation}
\mathcal{H}^{\rm s}(F) \longleftrightarrow \widetilde{Q}(F) \otimes \widetilde{Q}(F) 
\quad \mbox{and} \quad 
\mathcal{H}^{\rm w}(F) \longleftrightarrow \widetilde{Q}(F) \otimes \widetilde{Q}(F)/ \langle \tau, \rho \rangle.
\end{equation} 
\end{thm} 
This theorem is proved as Theorems~\ref{thm:HQs} and \ref{thm:HQw}. 

When $F$ is an oriented surface-link, a $1$-handle attached to $F$ is called an {\it oriented $1$-handle} if the orientation is compatible with the orientation of $F$ (Definition~\ref{def:orientedhandle}).  
Let 
$\mathcal{H}^{\rm s}(F)^{\rm ori}$ (or $\mathcal{H}^{\rm w}(F)^{\rm ori}$)  denote the 
set of strong (or weak) equivalence classes of oriented $1$-handles attached to $F$.   

\begin{thm}\label{thm:mainB}
Let $F$ be an oriented surface-link and $Q(F)$ the knot quandle of $F$.  
There are bijections
\begin{equation}
\mathcal{H}^{\rm s}(F)^{\rm ori} \longleftrightarrow Q(F) \otimes Q(F) 
\quad \mbox{and} \quad 
\mathcal{H}^{\rm w}(F)^{\rm ori} \longleftrightarrow Q(F) \otimes Q(F)/ \langle \tau \rangle. 
\end{equation} 
\end{thm}
This theorem is proved as Theorems~\ref{thm:HQsori} and \ref{thm:HQwori}. 

Theorems~\ref{thm:mainA} and \ref{thm:mainB} give complete classifications of all $1$-handles and oriented $1$-handles in terms of the tensor products of quandles.  
When $\widetilde{Q}(F)$ is finite, the tensor product $\widetilde{Q}(F) \otimes \widetilde{Q}(F)$ and its quotient $\widetilde{Q}(F) \otimes \widetilde{Q}(F)/ \langle \tau, \rho \rangle$ are finite and we can list the elements theoretically.  However, when  $\widetilde{Q}(F)$ is infinite or has a large cardinality, it is difficult to compute or list the elements of the tensor product and its quotient.  In such a case, it is useful to construct an invariant.  Theorems~\ref{thm:mainA} and \ref{thm:mainB} can be used to construct invariants of $1$-handles as follows:  

\begin{thm}\label{thm:mainC}
Let $F$ be a surface-link and $(\widetilde{Q}(F), \rho)$ the knot symmetric quandle of $F$.  
Let $f: (\widetilde{Q}(F), \rho) \to (X, \rho)$ be a symmetric quandle homomorphism to a finite symmetric quandle $(X, \rho)$.  
We have invariants of $1$-handles, 
\begin{equation}
\mathcal{H}^{\rm s}(F) \longrightarrow X \otimes X
\quad \mbox{and} \quad 
\mathcal{H}^{\rm w}(F) \longrightarrow X \otimes X/ \langle \tau, \rho \rangle.
\end{equation} 
\end{thm} 

\begin{thm}\label{thm:mainD}
Let $F$ be an oriented surface-link and $Q(F)$ the knot quandle of $F$.  
Let $f: Q(F) \to X$ be a quandle homomorphism to a finite quandle $X$.  
We have invariants of $1$-handles, 
\begin{equation}
\mathcal{H}^{\rm s}(F)^{\rm ori} \longrightarrow X \otimes X 
\quad \mbox{and} \quad 
\mathcal{H}^{\rm w}(F)^{\rm ori} \longrightarrow X \otimes X/ \langle \tau \rangle. 
\end{equation} 
\end{thm}

This paper is organized as follow: In Section~\ref{sect:quandles} we recall quandles and symmetric quandles. In Section~\ref{sect:tensor} the tensor product of a quandle is defined.  For a quandle $X$,  we give a method of computing the tensor product $D(X) \otimes D(X)$ of the symmetric double $D(X)$ 
from the tensor product $X \otimes X$  (Theorem~\ref{thm:tensordouble}).  
In Section~\ref{sect:involution} for a quandle $X$ or a symmetric quandle $(X, \rho)$, involutions $\tau$ and $\rho$ on the tensor product $X \otimes X$ and the quotients $X \otimes X/ \langle \tau \rangle$, $X \otimes X/ \langle \rho \rangle$ and $X \otimes X/ \langle \tau, \rho \rangle$ are defined.  
In Section~\ref{sect:dihedral} the tensor product $R_n \otimes R_n$ of the dihedral quandle $R_n$ is computed.  
In Section~\ref{sect:d:dihedral} the  tensor product $D(R_n) \otimes D(R_n)$ and its quotients for the symmetric double $D(R_n)$ of $R_n$ are computed.  
 In Section~\ref{sect:handle}  the definition of strong/weak equivalence of $1$-handles is given.  
We recall the notion of chords attached to a surface-link and a relationship between $1$-handles and chords.  
   In Section~\ref{sect:strong} we discuss strong equivalence classes of $1$-handles and show  a bijection from the strong equivalence classes of $1$-handles  (or oriented $1$-handles)  attached to $F$ to the tensor product 
$\widetilde{Q}(F) \otimes \widetilde{Q}(F)$ (or $Q(F) \otimes Q(F) $) (Theorems~\ref{thm:HQs} and \ref{thm:HQsori}).   
   In Section~\ref{sect:weak} we discuss weak equivalence classes of $1$-handles and show  a bijection from the weak equivalence classes of $1$-handles  (or oriented $1$-handles)  attached to $F$ to the quotient of the tensor product 
$\widetilde{Q}(F) \otimes \widetilde{Q}(F)/ \langle \tau, \rho \rangle$ (or $Q(F) \otimes Q(F)/ \langle \tau \rangle$) (Theorems~\ref{thm:HQw} and \ref{thm:HQwori}).   
   In Section~\ref{sect:example} we give examples and discuss invariants of $1$-handles.  

This work was supported by JSPS KAKENHI Grant Numbers JP19H01788 and JP17H06128.


\section{Quandles and symmetric quandles}
\label{sect:quandles}

A {\it quandle} is a set $X$ with two binary operations $\ast$ and $\overline{\ast}$ 
satisfying  the following three conditions: 
\begin{itemize} 
\item[(Q1)] For any $x \in X$, $x \ast x=x$.  
\item[(Q2)] For any $x, y \in X$, $(x \ast y)  \,\overline{\ast}\, y =x$ and 
$(x \,\overline{\ast}\, y)  \ast y =x$. 
\item[(Q3)] For any $x, y, z \in X$, $(x \ast y) \ast z = (x \ast z) \ast (y \ast z)$. 
\end{itemize} 
This notion was introduced independently by Joyce \cite{Joyce82} 
and Matveev \cite{Matveev82}. The three conditions correspond to three basic moves on knot diagrams, called Reidemeister moves.  
A {\it symmetric quandle} is a pair $(X, \rho)$ of a quandle $X$ and an involution $\rho: X \to X$, 
called a {\it good involution} on $X$,  satisfying the following conditions: 
\begin{itemize} 
\item[(S1)] For any $x, y \in X$,  $\rho(x \ast y) = \rho(x) \ast y$,  
\item[(S2)] For any $x, y \in X$,  $x \ast {\rho(y)} = x \, \overline{\ast}\, y $. 
\end{itemize} 
A {\it symmetric quandle homomorphism} $f : (X, \rho) \to (Y, \rho)$ means a quandle homomorphism 
$f : X \to Y$ with $f \circ \rho = \rho \circ f$. 
Refer to \cite{Kamada07, Kamada14A, KO10} for details on symmetric quandles and related topics. 

We here give some examples of quandles and symmetric quandles.  

\begin{example}{\rm 
(Dihedral quandles) 
Let $X= \Z / n \Z= \{ 0, 1, \cdots, n-1\}$.  
Define $a \ast b = a \,\overline{\ast}\, b = 2b -a$.  
Then $X = (X, \ast,  \,\overline{\ast}\, )$ is a quandle, which is called the 
{\it dihedral quandle} of order $n$.  It is denoted by $R_n$ in this paper.  

Let $\rho: R_n \to R_n$ be the identity map. 
It is a good involution, and we have a symmetric quandle $(R_n, \rho)$. 

Suppose that $n= 2m$ and let $\rho' : R_n \to R_n$ be the map sending $x$ to 
$ x + m$ for all $x \in R_n$. Then $\rho'$ is a good involution of $R_n$, and we have another symmetric quandle $(R_n, \rho')$ whose underlying quandle is the dihedral quandle.   Good involutions on dihedral quandles are classified in \cite{KO10}.   
}\end{example}

\begin{example}{\rm 
(Conjugation quandles) 
Let $G$ be a group.  
Let $X$ be $G$ as a set. 
Define $a \ast b = b^{-1} a b$ and $a \,\overline{\ast}\, b = b a b^{-1}$.  
Then $X = (X, \ast,  \,\overline{\ast}\, )$ is a quandle, which is called the 
{\it conjugation quandle} of $G$ (cf. \cite{Joyce82, Matveev82}).  

Let $\rho: X \to X$ be the map sending $x$ to $x^{-1}$ for all $x \in X$. It is a good involution on the conjugation quandle, and we have a symmetric quandle $(X, \rho)$ (cf. \cite{Kamada07, Kamada14A, KO10}).  
}\end{example}

\begin{example}\label{doublecover}{\rm 
{\rm (Symmetric doubles, cf. \cite{Kamada07, Kamada14A, KO10})
Let $X = (X, \ast,  \,\overline{\ast}\, )$ be a quandle. 
Let $X^+$ and $X^-$ be copies of $X$. 
For each $x \in X$, we denote by $x^+$ and $x^-$ the elements of $X^+$ and $X^-$ corresponding to $x$.  
Let $D(X) = X^+ \amalg X^-$ and  define  binary operations $\ast$ and $\overline{\ast}$ on 
$D(X)$ by 
\begin{equation}
\begin{array}{ll} 
x^\pm \ast y^+  = (x \ast y)^\pm,   \quad & 
x^\pm \ast y^- = (x  \,\overline{\ast}\,  y )^\pm,  \\ 
x^\pm \,\overline{\ast}\, y^+  = (x \,\overline{\ast}\, y)^\pm,  \quad &  
x^\pm \,\overline{\ast}\, y^- = (x  \ast  y )^\pm,  
\end{array}
\end{equation}
where $x^\pm$ means $x^+$ or $x^-$, respectively.  
Then $D(X) = (D(X), \ast,  \,\overline{\ast}\, )$ is a quandle, which we call the {\it symmetric double} of $X$.  
(The subquandle $X^+$ of $D(X)$ can be identified with the original quandle $X$  by identifying $x^+$ with $x$. Elements $x^+$ and $x^-$ are often denoted by $x$ and $\overline{x}$, respectively.)
The involution $\rho: D(X) \to D(X)$ interchanging  $x^+$ and  $x^-$ for all $x\in X$ is a good involution.  The symmetric quandle $(D(X), \rho)$, with this particular $\rho$, is also referred to as the {\it symmetric double} of $X$. 
}
}\end{example}

\begin{remark}\label{lem:symdouble} 
Let $X$ and $Y$ be quandles and let $f: X \to Y$ be a quandle homomorphism. 
Let $f$ also denote the map $f : D(X) \to D(Y)$ sending $x^\epsilon$ to $f(x)^\epsilon$ for $x \in X$ and $\epsilon \in \{+,-\}$.  Then $f : D(X) \to D(Y)$ is a quandle homomorphism.  Moreover, 
when $(D(X), \rho)$ and $(D(Y), \rho)$ are the symmetric doubles of $X$ and $Y$, 
the map $f : D(X) \to D(Y)$ is a symmetric quandle homomorphism, i.e., it is a quandle homomorphism with $f (\rho(a) ) = \rho (f(a))$ for $a \in D(X)$.  
\end{remark}

\begin{figure}[h]
\begin{center}
\mbox{\epsfxsize=3.2in \epsfbox{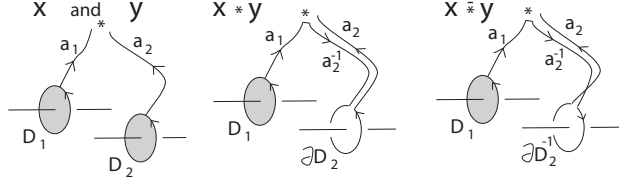}}
\end{center}
\vspace{-0.5cm}
\caption{Knot quandle operations}
\label{fgda.eps}
\end{figure}

\begin{example}\label{example:knotqdle}
{\rm (Knot quandles, cf. \cite{FennRourke92, Joyce82, Kamada14A, Matveev82})
Let $ K$ be an oriented closed $n$-manifold embedded in $\R^{n+2}$.  
Let ${Q}(K) = \{ [(D,a)] \}$ be the set of homotopy classes of 
pairs $(D, a)$ such that $D$ is 
a positively oriented meridian disk  of $K$ and $a$ is an arc 
in the knot exterior connecting $D$ and the base point. 

Define 
\begin{equation}
\begin{array}{lll} 
~[(D_1, a_1)] \ast [(D_2, a_2)] & = & [(D_1, a_1 \cdot a_2^{-1} \cdot \partial D_2 \cdot a_2)],  \\ 
~[(D_1, a_1)] \, \overline{\ast} \, ~  [(D_2, a_2)] & = & [(D_1, a_1 \cdot a_2^{-1} \cdot \partial (-D_2) \cdot a_2)]. 
\end{array}
\end{equation}

See Figure~\ref{fgda.eps}.  The {\it knot quandle } ${Q}(K)$ of $K$ is $({Q}(K), \ast, \overline{\ast})$. 
The knot quandle of $K$  can be computed from a diagram of $K$ (cf. \cite{CKS, Joyce82, Kamada14A, Matveev82}). 
}\end{example}

\begin{figure}[h]
\begin{center}
\mbox{\epsfxsize=2.5in \epsfbox{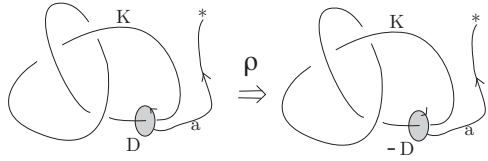}}
\end{center}
\vspace{-0.5cm}
\caption{The standard involution $\rho$}
\label{fgc.eps}
\end{figure}

\begin{example}\label{example:knotsymqdle}
{\rm (Full knot quandles and knot symmetric quandles, cf. \cite{Kamada07, Kamada14A, KO10})
Let $ K$ be a closed $n$-manifold embedded in $\R^{n+2}$, which may be non-orientable.  
Let $\widetilde{Q}(K) = \{ [(D,a)] \}$ be the set of homotopy classes of 
pairs $(D, a)$ such that $D$ is 
an oriented meridian disk of $K$ and $a$ is an arc 
in the knot exterior connecting $D$ and the base point. 
     Define binary operations $\ast$ and $\, \overline{\ast} \,$ on 
$\widetilde{Q}(K)$ by the same way as the operations of the knot quandle 
above.  Then $(\widetilde{Q}(K), \ast, \overline{\ast})$ is a quandle, which we call the  
 {\it full knot quandle} $\widetilde{Q}(K)$ of $K$. 

The {\it standard involution} is $\rho : \widetilde{Q}(K) \to \widetilde{Q}(K)$, 
$\rho [(D, a)] = [(-D, a)]$, which is  a good involution.  See Figure~\ref{fgc.eps}. 
The {\it knot symmetric quandle} of $K$ is the pair $(\widetilde{Q}(K), \rho)$ 
of the full knot quandle and its standard involution.   
}\end{example}

\begin{remark}\label{remark:double}{\rm 
When $K$ is oriented, the knot quandle $Q(K)$ is a subquandle of 
$\widetilde{Q}(K)$, and the knot symmetric quandle $(\widetilde{Q}(K), \rho)$  is identified with the symmetric double $(D(Q(K)), \rho)$ of $Q(K)$ (cf.  \cite{Kamada07, Kamada14A, KO10}).  
}\end{remark}

\vspace{0.0cm}

\begin{remark}\label{remark:group}{\rm 
For a quandle $X$,   the {\it associated group} of $X$ is defined by 
\begin{equation}
 As(X) = \langle x \in X \mid  x \ast y = y^{-1} x y \quad 
\mbox{for $x, y \in X$} 
\rangle. 
\end{equation}

It is known that for 
 an oriented closed $n$-manifold $K$ embedded in $\R^{n+2}$,   
the knot group $G(K)= \pi_1(\R^{n+2} \setminus K)$ is isomorphic to $As(Q(K))$ 
(cf. \cite{FennRourke92, Joyce82, Matveev82}). 

For a symmetric quandle $(X, \rho)$, the {\it associated group} of $(X, \rho)$ is 
defined by 
\begin{equation}
As(X, \rho) = \langle x \in X \mid  x \ast y = y^{-1} x y , ~ \rho(x)= x^{-1} \quad 
\mbox{for $x, y \in X$} 
\rangle. 
\end{equation} 

For a closed $n$-manifold $K$ embedded in $\R^{n+2}$, the  
knot group $G(K)$ is isomorphic to $As(\widetilde{Q}(K), \rho)$ 
(Theorem~4.5 of \cite{Kamada14A}). 
}\end{remark}


\section{Tensor products of quandles}
\label{sect:tensor}


Suppose that a group $G$ acts on a quandle $X_1$ from the right and  on a quandle $X_2$ from the left;   $x_1 \cdot (g g') = (x_1 \cdot g) \cdot g'$ and $(g g') \cdot x_2 = g \cdot ( g' \cdot x_2)$. 
Consider the equivalence relation on $X_1 \times X_2$ generated by 
 \begin{equation}
(x_1 \cdot g, \, x_2) \sim (x_1, \, g \cdot x_2) \quad \mbox{for $x_1 \in X_1, \, x_2 \in X_2,$ and $g \in G.$}
\end{equation}
\begin{defn}{\rm 
The {\it tensor product} of $X_1$ and $X_2$ over $G$, denoted by 
 $X_1 \otimes_G X_2$, is the set of all equivalence classes $[(x_1, x_2)]$ 
for $(x_1, x_2) \in X_1 \times X_2$.   An element $[(x_1, x_2)]$ of $X_1 \otimes_G X_2$ is also denoted by $[x_1, x_2]$,  $x_1 \otimes x_2$ or $(x_1, x_2)$. 
}\end{defn}

Let $X$ be a quandle.  Let $F(X)$ be the free group generated by $X$.  
The free group $F(X)$ acts on $X$ from the right by 
\begin{equation}
x \cdot g := (((x \ast^{\epsilon_1} x_1) \ast^{\epsilon_2} x_2) \cdots ) \ast^{\epsilon_m} x_m 
\end{equation}
for $x \in X$ and $g= x_1^{\epsilon_1} x_2^{\epsilon_2}\cdots x_m^{\epsilon_m} \in G$, 
where $\ast^1 = \ast$ and $\ast^{-1}= \overline{\ast}$.  
For example, $x \cdot ab c^{-1} d=  (((x \ast a) \ast b) \, \overline{\ast}\, c) \ast d$.  
This action is called the {\it canonical action} of $F(X)$ on $X$. 
The free group $F(X)$ acts on $X$ from the left by 
\begin{equation}
g \cdot x   :=  x \cdot g^{-1}   = 
(((x \ast^{-\epsilon_m} x_m) \ast^{-\epsilon_{m-1}} x_{m-1}) \cdots ) \ast^{-\epsilon_1} x_1  
\end{equation} 
for $x \in X$ and $g= x_1^{\epsilon_1} x_2^{\epsilon_2}\cdots x_m^{\epsilon_m} \in G$.

\begin{defn}{\rm 
The {\it canonical tensor product} of $X$ is the tensor product 
$X \otimes_{F(X)} X$ of $X$ over the canonical action of $F(X)$.   In this paper, 
we refer to the canonical tensor product as 
the {\it tensor product} of $X$ and denote it simply by $X \otimes X$.   
}\end{defn}

Note that 
\begin{equation}
\begin{array}{rl} 
X \otimes X & = X \times X /  (x_1 \cdot g, \, x_2) \sim (x_1, \, g \cdot x_2) \\
& = X \times X /  (x_1 \cdot g, \, g^{-1} \cdot x_2 ) \sim (x_1, \, x_2) \\
& = X \times X /  (x_1 \cdot g, \, x_2 \cdot g ) \sim (x_1, \, x_2). \\
\end{array}
\end{equation}
Thus, the tensor product $X \otimes X$ is the orbit 
space of $X \times X$ by the right action of $F(X)$ with  
$(x_1, \, x_2) \cdot g := (x_1 \cdot g, \, x_2 \cdot g)$.  
We also refer to this action of $F(X)$ on $X \times X$, $(x_1, \, x_2) \cdot g := (x_1 \cdot g, \, x_2 \cdot g)$, as the {\it canonical right action} of $F(X)$ on $X \times X$.

\begin{example}\label{example:R3}{\rm 
Let $R_3$ be the dihedral quandle of order $3$.  
As a set, $R_3 = \Z / 3 \Z = \{ 0,1,2\}   $ and $R_3 \times R_3 = \{ (i, j ) \mid i, j \in \{ 0,1,2\}  \}$.  
Note that $(i, j) \sim (i \cdot k, j \cdot k) =  
(i \ast k, j \ast k) =(2k -i, 2k - j)$ for any $i,j,k$.  

For example, 
$ (0, 0) \sim   (0\ast 1, 0\ast 1) = (2, 2) $   and  $ (0, 1) \sim  (0\ast 1, 1\ast 1) = (2, 1)$.  

Considering the canonical right actions of $0, 1$ and $2$ on $R_3 \times R_3$, we see that 
\begin{equation}
\begin{array}{ll} 
~[0,0] & = \{(0,0), (1,1), (2,2) \}   \\ 
~[0,1] & = \{(0,1), (0,2),  (1,0),  (1,2), (2,0), (2,1)\},  
 \end{array} 
\end{equation}
and 
\begin{equation}
R_3 \otimes R_3 =\{ [0,0], [0,1] \}. 
\end{equation}
}\end{example}

In Section~\ref{sect:dihedral} we compute the tensor product $R_n \otimes R_n$ for every dihedral qaundle $R_n$.

We introduce a method of computing $D(X) \otimes D(X)$ from $X \otimes X$.  

Let $X$ be a quandle and let $D(X)$ be 
the symmetric double.  Note that 
$D(X) =  X^+ \amalg X^- = \{ x^+, x^- \mid x \in X\}$.  
For each element $(x, y) \in X \times X$ and $\epsilon, \delta \in \{+, -\}$, let 
$(x, y)^{\epsilon, \delta}$ denote the element $(x^\epsilon, y^\delta) \in D(X) \times D(X)$.  
Note that every element of $D(X) \times D(X)$ can be written in this form. 
For a subset $E$ of $X \times X$ and $\epsilon, \delta \in \{+, -\}$, let 
$E^{\epsilon, \delta}$ denote the subset of $D(X) \times D(X)$ 
with 
\begin{equation}
E^{\epsilon, \delta} = \{  (x, y)^{\epsilon, \delta} \mid (x, y) \in E \}.  
\end{equation}

\begin{thm}\label{thm:tensordouble} 
For any quandle $X$, 
\begin{equation}
D(X) \otimes D(X) = \{ E^{\epsilon, \delta} \mid  E \in X \otimes X ~\mbox{\textrm{and}}~ 
\epsilon, \delta \in \{+,-\} \}.
\end{equation}
\end{thm}

We provide two lemmas.  

Let $\alpha^+ : F(X) \to F(D(X))$ be the homomorphism determined by 
$\alpha^+ (x) = x^+$ for $x \in X$.  (Then $\alpha^+ (x^{-1}) = (x^+)^{-1}$ for $x \in X$.)  

Let $\beta: F(D(X)) \to F(X)$ be the homomorphism determined by 
$\beta(x^+) = x$ and $\beta(x^-)= x^{-1}$ for $x \in X$.  
(Then $\beta((x^+)^{-1}) = x^{-1}$ and $\beta((x^-)^{-1}) = x$ for $x \in X$.)

\begin{lem}\label{lem:tdA} 
Let $x$, $x_0 \in X$,   $g \in F(X)$, $h \in F(D(X))$, and $\epsilon \in \{+,-\}$.  
\begin{itemize}
\item[(1)] 
If $x = x_0 \cdot g$, then  $x^\epsilon = (x_0)^\epsilon \cdot \alpha^+(g)$.  
\item[(2)] 
If $x^\epsilon = (x_0)^\epsilon \cdot h$, then   $x = x_0 \cdot \beta(h)$.  
\end{itemize}
\end{lem}

{\it Proof.}~ 
It is directly seen from the definition of the canonical action. \qed 

\begin{lem}\label{lem:tdB} 
Let $E \in X \otimes X$ and $E = [(x_0, y_0)]$.  Then 
$E^{\epsilon, \delta} = [(x_0, y_0)^{\epsilon, \delta}]$. 
\end{lem}

{\it Proof.}~  
Let $(x, y)^{\epsilon, \delta}$ be any element of  $E^{\epsilon, \delta}$.  
By definition, $(x, y) \in E$.  There is an element $g \in F(X)$ with $(x, y) 
= (x_0, y_0) \cdot g$.  By Lemma~\ref{lem:tdA}, $(x, y)^{\epsilon, \delta} 
= (x_0, y_0)^{\epsilon, \delta} \cdot \alpha_+(g)$.  Thus, 
$(x, y)^{\epsilon, \delta}  \in [(x_0, y_0)^{\epsilon, \delta}]$.  
Hence, $E^{\epsilon, \delta} \subset [(x_0, y_0)^{\epsilon, \delta}]$.  

Conversely, let $(x, y)^{\epsilon, \delta} \in D(X) \times D(X)$ 
be any element of $[(x_0, y_0)^{\epsilon, \delta}]$. 
There is an element $h \in F(D(X))$ with $(x, y)^{\epsilon, \delta} = 
(x_0, y_0)^{\epsilon, \delta} \cdot h$.  
By Lemma~\ref{lem:tdA},  $(x, y) = 
(x_0, y_0)  \cdot \beta(h) \in E$.  Thus, 
$(x, y)^{\epsilon, \delta}  \in E^{\epsilon, \delta}$. 
Hence,  $[(x_0, y_0)^{\epsilon, \delta}] \subset E^{\epsilon, \delta}$. \qed 

\vspace{2mm}
{\it Proof of Theorem~\ref{thm:tensordouble}.}~ 
It follows from Lemma~\ref{lem:tdB}. \qed 

\vspace{2mm}
The following is seen from Example~\ref{example:R3} and Theorem~\ref{thm:tensordouble}. 

\begin{example}\label{example:DR3}{\rm 
Let $D(R_3)$ be the symmetric double of $R_3$.  
Then 
\begin{equation}
D(R_3) \otimes D(R_3) =\{ 
[0, 0]^{\epsilon, \delta},  
[0, 1]^{\epsilon, \delta} \mid   \epsilon, \delta \in \{+,-\}    \}, 
\end{equation}
where 
\begin{equation}
\begin{array}{ll}
~[0, 0]^{\epsilon, \delta} & = \{ (0, 0)^{\epsilon, \delta}, (1, 1)^{\epsilon, \delta}, (2, 2)^{\epsilon, \delta} \},  \\  
~[0, 1]^{\epsilon, \delta} & = \{(0, 1)^{\epsilon, \delta}, (0, 2)^{\epsilon, \delta},  (1, 0)^{\epsilon, \delta},  (1, 2)^{\epsilon, \delta}, (2, 0)^{\epsilon, \delta}, (2, 1)^{\epsilon, \delta} \},  \quad  \\ 
\end{array} 
\end{equation}
for $\epsilon, \delta \in \{+,-\}$.  
Thus $D(R_3) \otimes D(R_3)$ consists of $8$ elements. 
}\end{example}


\section{Involutions on tensor products}
\label{sect:involution} 


In this section we introduce involutions on tensor products.  

Let $X \otimes X$ be the tensor product of a quandle $X$.  
Let $\tau: X \times X \to X \times X$ be the involution of $X \times X$ 
switching the components, $\tau(x_1, x_2)= (x_2, x_1)$.  Since $\tau$ commutes with 
the canonical right action of $F(X)$ on $X \times X$, it induces an involution on $X \otimes X$, 
which we denote also by 
\begin{equation}
\tau : X \otimes X \to X \otimes X, \quad [x_1, x_2] \mapsto [x_2, x_1]. 
\end{equation}

Let $\rho: X \to X$ be an involution.  
It induces an invitation $\rho ~ (= \rho \otimes \rho): X \times X \to X \times X$ with 
$\rho(x_1, x_2) = (\rho(x_1), \rho(x_2))$.  Note that $\rho: X \times X \to X \times X$ commutes with the canonical right action of $F(X)$ on $X \times X$ if and only if $\rho: X \to X$ satisfies the first condition (S1) of a good involution.  

Let $\rho: X \to X$ be a good involution. Then the involution $\rho: X \times X \to X \times X$ 
induces an  involution on $X \otimes X$, 
which we denote also by 
\begin{equation}
\rho : X \otimes X \to X \otimes X, \quad [x_1, x_2] \mapsto [\rho(x_1), \rho(x_2)]. 
\end{equation}

We  denote by $X \otimes X/ \langle \tau \rangle$ ($X \otimes X/ \langle \rho \rangle$ or $X \otimes X/ 
\langle \tau, \rho \rangle$, resp.)
the set of equivalence classes of $X \otimes X$ by the equivalence relation generated by $\tau$ 
(by $\rho$ or by $\tau$ and $\rho$, resp.). 

\begin{example}\label{example:R3inv}{\rm 
Let $R_3$ be the dihedral quandle of order $3$.  Recall that 
$R_3 \otimes R_3 =\{ [0,0], [0,1] \}$ as in Example~\ref{example:R3}. 
Since $\tau( [0,0] )= [0,0]$ and $\tau ([0,1]) = [0,1]$, we have 
\begin{equation}
R_3 \otimes R_3 / \langle \tau \rangle  =\{ \{[0,0]\}, \{[0,1]\} \}.   
\end{equation}
}\end{example}

In Section~\ref{sect:dihedral} we compute $R_n \otimes R_n/ \langle \tau \rangle$ for every dihedral qaundle $R_n$.

\begin{example}\label{example:DR3inv}{\rm 
Let $(D(R_3), \rho)$ be the symmetric double of $R_3$. 
Recall that $
D(R_3) \otimes D(R_3) =\{ 
[0, 0]^{\epsilon, \delta},   [0, 1]^{\epsilon, \delta},  \mid  \epsilon, \delta \in \{+,-\}  \}$ 
as in Example~\ref{example:DR3}. 
Note that for $k =0, 1$, 
\begin{equation}
\begin{array}{ll} 
\tau([0, k]^{+,+}) = [0, k]^{+,+}, \quad & \tau([0, k]^{-,-})= [0, k]^{-,-}, \\ 
\tau( [0, k]^{+,-}) =  [0, k]^{-,+},  &    \\ 
\rho([0, k]^{+,+}) = [0, k]^{-,-}, \quad & \rho( [0, k]^{+,-}) =  [0, k]^{-,+}.    \\
\end{array}
\end{equation}
Thus we see the following.  
\begin{itemize}
\item 
$D(R_3) \otimes D(R_3) / \langle \tau \rangle$ 
consists of 
$$ \{[0, k]^{+,+} \}, \{ [0, k]^{-,-} \}, \{ [0, k]^{+,-},  [0, k]^{-,+}\} \quad  \mbox{ for $ k=0,1 $. } $$ 
\item 
$D(R_3) \otimes D(R_3) / \langle \rho \rangle $ 
consists of 
$$ \{ [0, k]^{+,+}, [0, k]^{-,-} \},   \{ [0, k]^{+,-}, [0, k]^{-,+} \}  \quad \mbox{ for  $ k=0,1 $.} $$  
\item 
$D(R_3) \otimes D(R_3) / \langle \tau, \rho \rangle$ 
consists of 
$$ \{ [0, k]^{+,+}, [0, k]^{-,-} \},   \{ [0, k]^{+,-}, [0, k]^{-,+} \}  \quad \mbox{ for  $ k=0,1 $.} $$   
\end{itemize}
}\end{example}

The following lemma is straightforward and we omit the proof.    

\begin{lem}\label{lem:tensorhom} 
For a quandle homomorphism $f: X \to Y$, 
we have a map 
\begin{equation}
 f \otimes f : X \otimes X  \to Y \otimes Y, \quad [x_1, x_2] \mapsto [f(x_1), f(x_2)].  
\end{equation}
It induces a map $
X \otimes X / \langle \tau \rangle \to Y \otimes Y / \langle \tau \rangle
$.  Moreover, for a symmetric quandle homomorphism 
$f: (X, \rho) \to (Y, \rho)$, the map $  f \otimes f $ induces maps 
$X \otimes X / \langle \rho \rangle \to Y \otimes Y / \langle \rho \rangle$ and 
$X \otimes X / \langle \tau, \rho \rangle \to Y \otimes Y / \langle \tau, \rho \rangle$. 
\end{lem}


\section{Tensor products of dihedral quandles} \label{sect:dihedral} 


In this section we compute tensor products $R_n \otimes R_n$  
and the quotient set $R_n \otimes R_n/ \langle \tau \rangle$ for every dihedral quandle $R_n$.  

\begin{thm}\label{thm:dihedral:odd}
Let $R_n= R_{2m+1}$ be the dihedral quandle of order $n= 2m+1$.  
The tensor product $R_n \otimes R_n$ consists of $m+1$ elements, $
E(0),  \dots, E(m)$,  
where 
\begin{equation}
E(k) = \{ (i, i+k), (i, i-k) \mid i \in R_n \} \quad \mbox{ for $k =0, 1, \dots, m$.}\\ 
\end{equation}
\end{thm}

We prove this theorem later.

Since $\tau (E(k)) = E(k)$ for $k =0, 1, \dots, m$, we have the following. 

\begin{cor}\label{cor:dihedral:odd}
Let $R_n= R_{2m+1}$ be the dihedral quandle of order $n= 2m+1$.  
The quotient set $R_n \otimes R_n / \langle \tau \rangle$ consists of $m+1$ elements, 
$\{ E(k) \}$ for  $k =0, 1, \dots, m$. 
\end{cor}

When $n=3$ and $m=1$, 
$E(0) = [0,0]$ and $E(1)=[0,1]$ in Examples~\ref{example:R3} and ~\ref{example:R3inv}. 

Let $R_n= R_{2m}$ be the dihedral quandle of order $n= 2m$.  
We denote by $R_{n, {\rm even}}$ (or $R_{n, {\rm odd}}$) 
the subset of $R_n$ consisting of the elements represented by even numbers 
(or odd numbers).  For example, when $n=8$, then 
$R_{8, {\rm even}}=\{ 0, 2, 4, 6\}$ and $R_{8, {\rm odd}}=\{ 1, 3, 5, 7\}$.  

 \begin{thm}\label{thm:dihedral:even}
Let $R_n= R_{2m}$ be the dihedral quandle of order $n= 2m$.  
The tensor product $R_n \otimes R_n$ consists of $n+2$ ($=2m+2$) elements, 
\begin{equation}
E(k)_0 \quad \mbox{and} \quad  E(k)_1 \quad \mbox{for $k=0, 1, \dots, m$,}
\end{equation}
where 
\begin{equation}
\begin{array}{ccll} 
E(k)_0  &  = &  \{ (i,i+k), (i, i-k) \mid  i \in R_{n, {\rm even}}  \}  \quad & \mbox{for $k=0, 1, \dots, m$,} \\ 
E(k)_1 &  = &  \{ (i, i+k), (i, i-k) \mid i \in R_{n, {\rm odd}} \}  \quad & \mbox{for $k=0, 1, \dots, m$.} \\ 
\end{array}
\end{equation}
\end{thm}

We prove this theorem later. 

It is obvious that 
\begin{equation}
\begin{array}{ll} 
\tau (E(k)_0) = E(k)_0, \quad \tau (E(k)_1) = E(k)_1 \quad & \mbox{for $k=0, 2, 4, \dots \leq m$,} \\ 
\tau (E(k)_0) = E(k)_1     & \mbox{for $k=1, 3, 5, \dots \leq m$.} \\ 
\end{array} 
\end{equation}

\begin{cor}\label{cor:dihedral:even}
Let $R_n= R_{2m}$ be the dihedral quandle of order $n= 2m$.  
The quotient set $R_n \otimes R_n / \langle \tau \rangle$ consists of $3(m+1)/2$ elements (if $m$ is odd) or 
$3m/2 + 2$ elements (if $m$ is even), which are listed below:  
\begin{equation}
\begin{array}{ll} 
\{ E(k)_0\}, \quad \{ E(k)_1\} \quad & \mbox{for $k= 0, 2, 4, \dots \leq m$}, \\ 
\{ E(k)_0, E(k)_1\}                 & \mbox{for $k= 1, 3, 5, \dots \leq m$}. \\ 
\end{array}
\end{equation}
\end{cor}

Let $d_n: R_n \times R_n \to \N_0 = \N \cup \{ 0 \}$ be the distance on $R_n = \Z / n \Z$ induced from the usual distance $d_0 : \Z \times \Z \to \N_0$ on $\Z$ with $d_0(x_1, x_2) = | x_1 - x_2 |$.  Since 
$d_n(x_1 \ast x, x_2 \ast x) = d_n(x_1 \,\overline{\ast}\, x, x_2 \,\overline{\ast}\, x) = d(x_1, x_2)$ for any $x_1, x_2, x \in R_n$, we have 
\begin{equation}\label{eqn:distance}
d_n(x_1 \cdot g, x_2 \cdot g) = d_n(x_1, x_2) \quad \mbox{ for any $x_1, x_2 \in R_n$ and $g \in F(R_n)$.} 
\end{equation}

\begin{lem}\label{lem:odd} 
Let $R_n= R_{2m+1}$ be the dihedral quandle of order $n= 2m+1$.  
Let $ k \in \{ 0, 1, \dots, m \}$. 
If $d_n(x, y) = k$, then there exists an element $g \in F(R_n)$ with $x \cdot g = 0$ and 
$y \cdot g = k$. 
\end{lem} 

{\it Proof.}~ 
It is known that $R_n$ with $n=2m+1$ is connected as a quandle, i.e., for any $x \in R_n$ there exists  an element $g \in F(R_n)$ with $x \cdot g = 0$.  In fact, let $g= x/2 \in R_n$ if $x$ is even, or $g= (n+x)/2 \in R_n$ if $x$ is odd. Then $x \ast g= 0$.  

Thus we have the assertion with $k=0$.  

Suppose $0 < k \leq m$ and $d_n(x, y) = k$.  Let $g_1 \in F(R_n)$ be an element with $x \cdot g_1 = 0$. By (\ref{eqn:distance}), we have $y \cdot g_1 = k$ or $-k$.  If $y \cdot g_1 = k$, then $g=g_1$ is a desired one. If $y \cdot g_1 = -k$, then $g= g_1 0$ is a desired one. \qed 

\vspace{2mm}
{\it Proof of Theorem~\ref{thm:dihedral:odd}.}~ 
For any $(x, y) \in R_n \times R_n$, $d_n(x, y) \in \{0, 1, \dots, m\}$.  Since 
\begin{equation}
E(k) = \{ (x, y) \in R_n \times R_n \mid d_n(x, y) =k \}, 
\end{equation}
we see that $E(0), E(1), \dots, E(m)$ cover $R_n \times R_n$ and 
$E(k) \cap E(k') = \emptyset$ for $k \neq k'$.  
By Lemma~\ref{lem:odd}, $E(k) \subset [0,k]$.  
Conversely, by (\ref{eqn:distance}), we see that $[0,k] \subset E(k)$. \qed 

\begin{lem}\label{lem:even} 
Let $R_n= R_{2m}$ be the dihedral quandle of order $n= 2m$.  
Let $ k \in \{ 0, 1, \dots, m \}$. 
If $d_n(x, y) = k$ and if $x$ is even (or odd), then there exists an element $g \in F(R_n)$ with $x \cdot g = 0$ and $y \cdot g = k$ (or $x \cdot g = 1$ and $y \cdot g = 1+ k$).  
\end{lem} 

{\it Proof.}~ 
It is known that $R_n$ with $n=2m$ has two connected components $R_{n, {\rm even}}$ and $R_{n, {\rm odd}}$, i.e., $x \in R_{n, {\rm even}}$ (or $x \in R_{n, {\rm odd}}$) if and only if  
 there exists  an element $g \in F(R_n)$ with $x \cdot g = 0$ (or $x \cdot g = 1$).  
In fact, if $x \in R_{n, {\rm even}}$ (or $x \in R_{n, {\rm odd}}$), then 
let $g= x/2$ (or $g= (1+x)/2$) and we have $x \cdot g = 0$ (or $x \cdot g = 1$).   
Conversely, if $x \cdot g = 0$ (or $x \cdot g = 1$) for some $g \in F(R_n)$, then 
 $x \in R_{n, {\rm even}}$ (or $x \in R_{n, {\rm odd}}$).  

Thus we have the assertion with $k=0$.  

Suppose $0 < k \leq m$, $d_n(x, y) = k$ and $x \in R_{n, {\rm even}}$ (or $x \in R_{n, {\rm odd}}$).  
Let $g_1 \in F(R_n)$ be an element with $x \cdot g_1 = 0$ (or $x \cdot g_1 = 1$).   
By (\ref{eqn:distance}), we have $y \cdot g_1 = k$ or $-k$ (or $y \cdot g_1 = 1+k$ or $1-k$).  
If $y \cdot g_1 = k$ (or $y \cdot g_1 = 1+k$), then $g=g_1$ is a desired one. If $y \cdot g_1 = -k$ (or $y \cdot g_1 = 1-k$), then $g= g_1 0$ (or $g= g_1 1$) is a desired one. \qed 

\vspace{2mm}
{\it Proof of Theorem~\ref{thm:dihedral:even}.}~ 
Since 
\begin{equation}
\begin{array}{ll} 
E(k)_0 & = \{ (x, y) \in R_n \times R_n \mid d_n(x, y) =k ~\mbox{and}~ x \in R_{n, {\rm even}} \} ~\mbox{and} \\ 
E(k)_1 & = \{ (x, y) \in R_n \times R_n \mid d_n(x, y) =k ~\mbox{and}~ x \in R_{n, {\rm odd}} \}, 
\end{array} 
\end{equation}
we see that $E(k)_0, E(k)_1$ for $k=0, 1, \cdots, m$  cover $R_n \times R_n$ and 
$E(k)_{i} \cap E(k')_{i'} = \emptyset$ unless $k= k'$ and $i = i'$.    
By Lemma~\ref{lem:odd}, $E(k)_0 \subset [0,k]$ and $E(k)_1 \subset [1,1+k]$.  
Conversely, by (\ref{eqn:distance}) and by the fact that $R_{n, {\rm even}}$ and $R_{n, {\rm odd}}$ are the connected components of $R_n$, 
we see that $[0,k] \subset E(k)_0$ and $[1,1+k] \subset E(k)_1$. \qed


\section{Tensor products of symmetric doubles of dihedral quandles} \label{sect:d:dihedral}


In this section we discuss tensor products of symmetric doubles of dihedral quandles, and their quotient sets by involutions.  

Let $R_n$ be the dihedral quandle of order $n$ and  $D(R_n)$  
its symmetric double.  Note that 
$D(R_n) = \{ i^+, i^- \mid i \in R_n\}$.  
We use the notation in Theorem~\ref{thm:tensordouble}.

\begin{thm}\label{thm:d:dihedral:odd}
Let $R_n= R_{2m+1}$ be the dihedral quandle of order $n= 2m+1$.  
The tensor product $D(R_n) \otimes D(R_n)$ consists of $4(m+1)$  ($=2n+2$) elements, 
$E(k)^{\epsilon, \delta}$ for $k=0,1, \dots, m$ and $\epsilon, \delta \in \{+, -\}$, 
where 
\begin{equation}
E(k)^{\epsilon, \delta} = \{ (i, i+k)^{\epsilon, \delta}, (i, i-k)^{\epsilon, \delta} \mid i \in R_n \} \quad \mbox{ for $k =0, 1, \dots, m$.}\\ 
\end{equation}
\end{thm}

{\it Proof.}~ 
It is a consequence of Theorems~\ref{thm:tensordouble} and 
\ref{thm:dihedral:odd}. \qed 

\vspace{2mm}
It is obvious from the definition of $E(k)^{\epsilon, \delta}$ that 
\begin{equation}
\begin{array}{l}
\tau (E(k)^{+,+}) = E(k)^{+,+}, \quad \tau (E(k)^{-,-}) = E(k)^{-,-}, \\ 
\tau (E(k)^{+,-}) = E(k)^{-,+},  \\ 
\rho (E(k)^{+,+}) = E(k)^{-,-}, \quad \rho (E(k)^{+,-}) = E(k)^{-,+},  
\end{array}
\end{equation}
for $k=0, 1, \dots, m$. 

Thus we have the following. 

\begin{cor}\label{cor:d:dihedral:odd:tau}
Let $R_n= R_{2m+1}$ be the dihedral quandle of order $n= 2m+1$.  
The quotient set $D(R_n) \otimes D(R_n) / \langle \tau \rangle$ consists of $3(m+1)$ elements, 
\begin{equation}
\{ E(k)^{+,+} \}, \{ E(k)^{-,-} \}, \{ E(k)^{+,-}, E(k)^{-,+} \} 
\quad \mbox{for $k=0, 1, \dots,  m$}. 
\end{equation}
\end{cor}

\begin{cor}\label{cor:d:dihedral:odd:rho}
Let $R_n= R_{2m+1}$ be the dihedral quandle of order $n= 2m+1$.  
The quotient set $D(R_n) \otimes D(R_n) / \langle \rho \rangle$ consists of $2(m+1)$  ($=n+1$) elements, 
\begin{equation}
\{ E(k)^{+,+}, E(k)^{-,-} \}, \{ E(k)^{+,-}, E(k)^{-,+} \} 
\quad \mbox{for $k=0, 1, \dots,  m$}. 
\end{equation}
\end{cor}

\begin{cor}\label{cor:d:dihedral:odd:taurho}
Let $R_n= R_{2m+1}$ be the dihedral quandle of order $n= 2m+1$.  
The quotient set $D(R_n) \otimes D(R_n) / \langle \tau, \rho \rangle$ consists of $2(m+1)$  ($=n+1$) elements, 
\begin{equation}
\{ E(k)^{+,+}, E(k)^{-,-} \}, \{ E(k)^{+,-}, E(k)^{-,+} \} 
\quad \mbox{for $k=0, 1, \dots,  m$}. 
\end{equation}
\end{cor}

When $n=3$ and $m=1$, 
$E(k)^{\epsilon, \delta} = [0, k]^{\epsilon, \delta}$ for $k=0, 1$ and $\epsilon, \delta \in \{+,-\}$ 
in Examples~\ref{example:DR3} and ~\ref{example:DR3inv}.

\begin{thm}\label{thm:d:dihedral:even}
Let $R_n= R_{2m}$ be the dihedral quandle of order $n= 2m$.  
The tensor product $D(R_n) \otimes D(R_n)$ consists of $8(m+1)$  ($=4n+8$) elements, 
$E(k)_0^{\epsilon, \delta}$ and $E(k)_1^{\epsilon, \delta}$ 
 for $k=0,1, \dots, m$ and $\epsilon, \delta \in \{+, -\}$, 
where 
\begin{equation}
\begin{array}{ccl} 
E(k)_0^{\epsilon, \delta}  &  = &  
\{ (i,i+k)^{\epsilon, \delta}, (i, i-k)^{\epsilon, \delta} \mid  i \in R_{n, {\rm even}} \},    \\ 
E(k)_1^{\epsilon, \delta} &  = &  
\{ (i, i+k)^{\epsilon, \delta}, (i, i-k)^{\epsilon, \delta} \mid i \in R_{n, {\rm odd}} \}   \\ 
\end{array}
\end{equation}
for $k=0, 1, \dots, m$.  
\end{thm}

{\it Proof.}~ 
It is a consequence of Theorems~\ref{thm:tensordouble} and 
\ref{thm:dihedral:even}. \qed 

\vspace{2mm}

The following corollaries are easily verified and we omit the proofs.  

\begin{cor}\label{cor:d:dihedral:even:tau}
Let $R_n= R_{2m}$ be the dihedral quandle of order $n= 2m$.  
The quotient set $D(R_n) \otimes D(R_n) / \langle \tau \rangle$ 
consists of $5(m+1)$ elements (if $m$ is odd) or 
$5m+6$ elements (if $m$ is even), which are listed below: 
\begin{equation}
\begin{array}{lll} 
\{ E(k)_0^{+,+}  \}, \quad & 
\{ E(k)_0^{-,-}  \}, \quad & 
\{ E(k)_0^{+,-},   E(k)_0^{-,+} \} \\ 
\{ E(k)_1^{+,+}  \}, \quad & 
\{ E(k)_1^{-,-}  \}, \quad & 
\{ E(k)_1^{+,-},   E(k)_1^{-,+} \} \\ 
\end{array} 
\end{equation}
for $k=0, 2, 4, \dots \leq m$, and 
\begin{equation}
\begin{array}{ll} 
\{ E(k)_0^{+,+}, E(k)_1^{+,+} \}, \quad & \{ E(k)_0^{-,-},  E(k)_1^{-,-} \}, \\ 
\{ E(k)_0^{+,-}, E(k)_1^{-,+} \},       & \{ E(k)_1^{+,-},  E(k)_0^{-,+} \}  \\ 
\end{array} 
\end{equation}
for $k=1, 3, 5, \dots \leq m$.  
\end{cor}

\begin{cor}\label{cor:d:dihedral:even:rho}
Let $R_n= R_{2m}$ be the dihedral quandle of order $n= 2m$.  
The quotient set $D(R_n) \otimes D(R_n) / \langle \rho \rangle$ 
consists of $4(m+1)$ ($=2n+4$) elements, which are listed below: 
\begin{equation}
\begin{array}{ll} 
\{ E(k)_0^{+,+}, E(k)_0^{-,-} \}, \quad & \{ E(k)_0^{+,-},  E(k)_0^{-,+} \}, \\ 
\{ E(k)_1^{+,+}, E(k)_1^{-,-} \},       & \{ E(k)_1^{+,-},  E(k)_1^{-,+} \}  \\ 
\end{array} 
\end{equation}
for $k=0, 1,  \dots,  m$.  
\end{cor}

\begin{cor}\label{cor:d:dihedral:even:taurho}
Let $R_n= R_{2m}$ be the dihedral quandle of order $n= 2m$.  
The quotient set $D(R_n) \otimes D(R_n) / \langle \tau, \rho \rangle$ 
consists of $3(m+1)$ elements (if $m$ is odd) or 
$3m+4$ elements (if $m$ is even), which are listed below: 
\begin{equation}
\begin{array}{ll} 
\{ E(k)_0^{+,+}, E(k)_0^{-,-} \}, \quad & \{ E(k)_0^{+,-},  E(k)_0^{-,+} \}, \\ 
\{ E(k)_1^{+,+}, E(k)_1^{-,-} \}, \quad & \{ E(k)_1^{+,-},  E(k)_1^{-,+} \} \\ 
\end{array} 
\end{equation}
for $k=0, 2, 4, \dots \leq m$, and 
\begin{equation}
\begin{array}{l} 
\{ E(k)_0^{+,+}, E(k)_0^{-,-},  E(k)_1^{+,+},  E(k)_1^{-,-} \}, \\ 
\{ E(k)_0^{+,-}, E(k)_0^{-,+},   E(k)_1^{+,-},  E(k)_1^{-,+} \}  \\ 
\end{array} 
\end{equation}
for $k=1, 3, 5, \dots \leq m$.  
\end{cor}


\section{$1$-handles and chords attached to surface-links}
\label{sect:handle} 


Let $F$ be a surface-link.  
A {\it $1$-handle} attached to $F$ means an embedding $h: [0,1] \times B^2 \to \R^4$ with $F \cap h([0,1] \times B^2) = h(\{0,1\} \times B^2)$.  The restriction of $h$ to $[0,1] \times \{0\}$ $(=[0,1])$  is denoted by $a_h: [0,1] \to \R^4$ and 
called the {\it core map}.  

\begin{defn}{\rm 
Two $1$-handles $h$ and $h'$ attached to $F$ are {\it strongly equivalent} 
(or {\it weakly equivalent}, resp.)  
if there is an ambient isotopy of $\R^4$ 
carrying $h$ to $h'$ as a map 
(or carrying the image $h([0,1] \times B^2)$ to $h'([0,1] \times B^2)$ as a subset of $\R^4$, resp.)  
and keeping $F$ setwise fixed.  
}\end{defn}

We denote by $\mathcal{H}^{\rm s}(F)$ (or $\mathcal{H}^{\rm w}(F)$)   
the set of strong (or weak) equivalence classes of $1$-handles attached to $F$.  

Let $F$ be a surface-link in $\R^4$.  

\begin{itemize}
\item 
A {\it chord} attached to $F$ is 
a path $a: [0,1] \to \R^4$ such that  $a(t) \in F$ if and only if  $t \in \{0,1\}$.  
\item 
An {\it initial root} (or a {\it terminal root})  
of $a$ is a regular neighborhood in $F$ of the initial point $a(0)$ (or the terminal point  $a(1)$).   
\item 
An {\it initial root orientation} (or a {\it terminal root orientation}) 
of $a$  is an orientation of an initial root (or a terminal root).  
A {\it root orientaion} of $a$ means a pair of an initial root orientation and a terminal root orientaion of $a$. 
\item 
A {\it root-oriented chord}  is a chord equipped with a root orientation or a pair $(a, o)$ of a chord $a$ and a root orientation $o$ of $a$. 
\end{itemize}

\begin{figure}[h]
\begin{center}
\mbox{\epsfxsize=3.8in \epsfbox{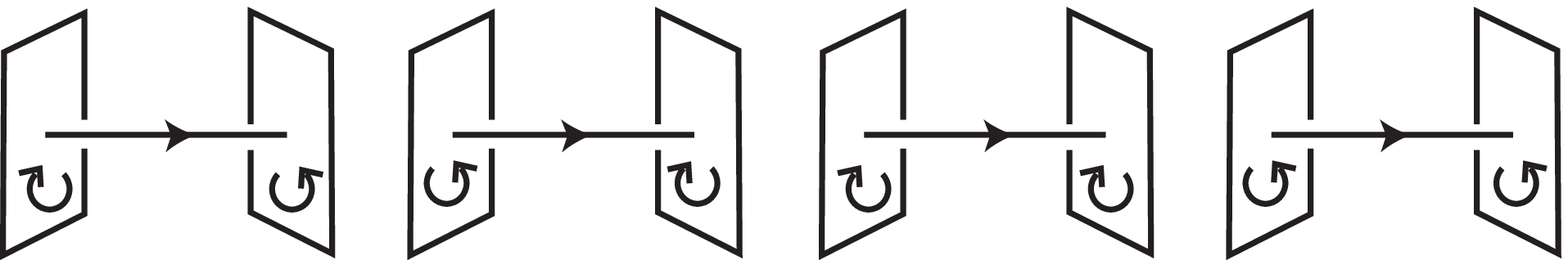}}
\put(-248,30){$a$}
\put(-177,30){$a$}
\put(-103,30){$a$}
\put(-33,30){$a$}
\put(-253,0){$(1)$}
\put(-182,0){$(2)$}
\put(-108,0){$(3)$}
\put(-38,0){$(4)$}
\end{center}
\vspace{0.0cm}
\caption{There are $4$ possible root orientations for a given chord $a$.}
\label{fgeb.eps}
\end{figure}

When a chord $a$ is given, there are $4$ possibilities for root orientations (Figure~\ref{fgeb.eps}). 

\begin{defn}
Two chords attached to $F$ are {\it strongly equivalent} 
or {\it homotopic} 
if they are homotopic as chords attached to $F$. 
Two root-oriented chords attached to $F$ are  {\it strongly equivalent} 
or {\it homotopic} 
if they are homotopic as root-oriented chords attached to $F$. 
\end{defn}

For a chord $a : [0,1] \to \R^4$ attached to $F$, we denote by 
${\rm rev}(a): [0,1]\to \R^4$ the chord defined by 
${\rm rev}(a) (t) = a(1-t)$.   For a root orientation $o$ of a chord $a$, 
we denote by $-o$ the reversed orientation of $o$. 

\begin{defn}
Two chords $a$ and $a'$ attached to $F$ are {\it weakly equivalent} 
if $a$ is homotopic to $a'$ or ${\rm rev}(a')$.  
Two root-oriented chords $(a, o)$ and $(a', o')$ attached to $F$  are {\it weakly equivalent} if $(a, o)$ is homotopic to $(a', o')$,  $(a', -o')$, 
$({\rm rev}(a'), o')$, or  $({\rm rev}(a'), -o')$.  
\end{defn}

Let $\mathcal{C}^{\rm s}(F)$ (or $\mathcal{C}^{\rm w}(F)$) denote the set of 
strong (or weak) equivalence classes 
of chords attached to $F$, 
and let $\mathcal{C}_{\rm root}^{\rm s}(F)$ (or $\mathcal{C}_{\rm root}^{\rm w}(F)$) denote 
the set of 
strong (or weak) equivalence classes  
of root-oriented chords attached to $F$.   

For a $1$-handle $h$ attached to a surface-link $F$, the core map $a_h$ is a chord attached to $F$.    
If $h$ and $h'$ are strongly (or weakly) equivalent, then $a_h$ and $a_{h'}$ are strongly (or weakly) equivalent.  
Thus we have maps  
\begin{equation}
\phi_0:  \mathcal{H}^{\rm s}(F) \to \mathcal{C}^{\rm s}(F) 
\quad \mbox{and} \quad 
\psi_0:  \mathcal{H}^{\rm w}(F) \to \mathcal{C}^{\rm w}(F) 
\end{equation}
sending $[h]$ to $[a_h]$.

Let $F$ be a surface-link, $h$  a $1$-handle attached to $F$ and 
$a_h$ the core map of $h$.  
Let $o_h$ be the root orientation of $a_h$ which is the reverse of the orientation 
induced from that of 
$\partial ([0,1] \times B^2)$ by $h$.  
We call $(a_h, o_h)$ the {\it root-oriented chord associated with } $h$. 

The following theorem is given  in \cite{Kamada14B}. 

\begin{thm}[\cite{Kamada14B}]\label{thm:HCsw} 
Two $1$-handles 
$h$ and $h'$ are strongly (or weakly) equivalent if and only if the  root-oriented chords $(a_h, o_h)$ and $(a_{h'}, o_{h'})$ associated with $h$ and $h'$ are strongly (or weakly) equivalent. 
In other words, the maps 
\begin{equation}
\phi:  \mathcal{H}^{\rm s}(F) \to \mathcal{C}_{\rm root}^{\rm s}(F)
\quad \mbox{and} \quad 
\psi:  \mathcal{H}^{\rm w}(F) \to \mathcal{C}_{\rm root}^{\rm w}(F)
\end{equation}
sending $[h]$ to $[(a_h, o_h)]$ are bijections.   
\end{thm}

In the rest of this section, 
let $F$ be an oriented surface-link, and we put a focus on oriented $1$-handles attached to $F$.  

\begin{defn}\label{def:orientedhandle}
A $1$-handle $h$ attached to $F$ is {\it oriented} if the orientation of 
$h (\{0, 1\} \times B^2)$ induced from the orientation of 
$\partial ([0,1] \times B^2)$ is opposite to that induced from the orientation of $F$.  
\end{defn}

In other words, when $(a_h, o_h)$ is the root-oriented chord associated with $h$, 
a $1$-handle $h$ attached to $F$ is oriented if and only if 
the root orientation $o_h$ matches the orientation of $F$.

Let  $\mathcal{H}^{\rm s}(F)^{\rm ori}$ 
(or $\mathcal{H}^{\rm w}(F)^{\rm ori}$) denote 
the subset of  $\mathcal{H}^{\rm s}(F)$ 
 (or $\mathcal{H}^{\rm w}(F)$) 
 consisting of the strong (or weak) equivalence classes of oriented $1$-handles attached to $F$.

The following is observed in \cite{Boyle, HK}, cf. \cite{Kamada14B}.  

\begin{thm}[\cite{Boyle, HK},  cf. \cite{Kamada14B}] \label{thm:A1} 
Let $F$ be an oriented surface-link. 
Two oriented $1$-handles $h$ and $h'$ attached to $F$ are strongly (or weakly) equivalent if and only if their chords $a_h$ and $a_{h'}$ are strongly (or weakly) equivalence.  In other words, 
the restriction maps 
of $\phi_0$ and $\psi_0$ to  $\mathcal{H}^{\rm s}(F)^{\rm ori}$ and 
$\mathcal{H}^{\rm w}(F)^{\rm ori}$,  
\begin{equation}
\phi_0^{\rm ori} :  \mathcal{H}^{\rm s}(F)^{\rm ori}  \to \mathcal{C}^{\rm s}(F)
\quad \mbox{and} \quad 
\psi_0^{\rm ori} :  \mathcal{H}^{\rm w}(F)^{\rm ori}  \to \mathcal{C}^{\rm w}(F),  
\end{equation}
sending $[h]$ to $[a_h]$, are  bijections.  
\end{thm}

\begin{remark}\label{remark:A}
For an oriented surface-link $F$,  
let $\mathcal{C}_{\rm root}^{\rm s}(F)^{\rm ori}$ (or $\mathcal{C}_{\rm root}^{\rm w}(F)^{\rm ori}$) 
denote the subset of 
$\mathcal{C}_{\rm root}^{\rm s}(F)$ (or $\mathcal{C}_{\rm root}^{\rm w}(F)$)
such that the root orientations match the orientation of $F$.  
There are natural bijections  
\begin{equation}
\mathcal{C}^{\rm s}(F) \to \mathcal{C}_{\rm root}^{\rm s}(F)^{\rm ori} 
\quad \mbox{and} \quad 
\mathcal{C}^{\rm w}(F) \to \mathcal{C}_{\rm root}^{\rm w}(F)^{\rm ori} 
\end{equation}
sending $[a]$ to $[(a, o_F)]$ where $o_F$ is the root orientation of $a$ which matches the orientation of $F$. 
Then the bijections $\phi_0^{\rm ori}$ and $\psi_0^{\rm ori}$ are regarded as the bijections 
\begin{equation}
\phi^{\rm ori} :  \mathcal{H}^{\rm s}(F)^{\rm ori}  \to \mathcal{C}_{\rm root}^{\rm s}(F)^{\rm ori}
\quad \mbox{and} \quad 
\psi^{\rm ori} :  \mathcal{H}^{\rm w}(F)^{\rm ori}  \to \mathcal{C}_{\rm root}^{\rm w}(F)^{\rm ori},   
\end{equation}
which are restrictions of the bijections $\phi$ and $\psi$ in Theorem~\ref{thm:HCsw}.  
\end{remark}


\section{$1$-handles and chords up to strong equivalence} 
\label{sect:strong}


In this section, we consider $1$-handles and chords up to strong equivalence.  
For a surface-link $F$, let $\widetilde{Q}(F)$ be the full knot quandle of $F$. 

Let $(D, a)$ be a pair of an oriented meridian disk $D$ of $F$ and a path 
$a:[0,1] \to \R^4$ in the knot exterior 
connecting $D$ and the base point, say $q$.  
Let $\tilde{a}: [0,1] \to \R^4$ be a path which is the concatenation of 
the path $a$ and 
a path 
in $D$ starting at the center of $D$ and terminating at $a(0)$. Then $\tilde{a}: [0,1] \to \R^4$ is a path in $\R^4$ such that 
$\tilde{a}(t) \in F$ if and only if $t=0$ and $\tilde{a}(1)= q$.  
We call $\tilde{a}$ an {\it extended path} of $a$.  
The initial point of $\tilde{a}$ is a point of $F$.  Using the orientation of $D$ and the standard orientation of $\R^4$, we can give an orientation of a regular neighborhood  in $F$ of the initial point of $\tilde{a}$, which we call the {\it root orientation} determined from the orientation of $D$ and we denote it by $o_{(D,a)}$. 

We  define a map 
$\mu: \widetilde{Q}(F) \otimes \widetilde{Q}(F) \to \mathcal{C}_{\rm root}^{\rm s}(F)$.  

Let $x_1 = [(D_1, a_1)]$ and $x_2 = [(D_2, a_2)]$ be elements of $\widetilde{Q}(F)$.  
Let $\tilde{a_1}$ and $\tilde{a_2}$ be extended paths of $a_1$ and $a_2$, 
and 
$o_{(D_1,a_1)}$ and $o_{(D_2,a_2)}$ be the root orientations determined from 
the orientations of $D_1$ and $D_2$.  
The concatenation of $\tilde{a_1}$ and $\tilde{a_2}^{-1}$ is 
a chord attached to $F$.  
We denote it by $a(x_1, x_2)$ and call it {\it the chord associated with 
$(x_1, x_2)$}.  We choose a root orientation of $a(x_1, x_2)$ such that the initial root orientation is $o_{(D_1,a_1)}$ and the terminal root orientation is $o_{(D_2,a_2)}$.  
Then  the strong equivalence class of the root-oriented chord is uniquely determined from 
$x_1 = [(D_1, a_1)]$ and $x_2 = [(D_2, a_2)]$.  
We denote this root-oriented chord by $(a(x_1, x_2), o(x_1, x_2))$ and call it the {\it root-oriented chord associated with } $(x_1, x_2)$.

\begin{thm}\label{thm:QCs} 
Let 
\begin{equation}
 \widetilde{Q}(F) \times \widetilde{Q}(F) \to \mathcal{C}_{\rm root}^{\rm s}(F) 
\end{equation}
be the map 
sending $(x_1, x_2) \in \widetilde{Q}(F) \times \widetilde{Q}(F)$ to 
the strong equivalence class of $(a(x_1, x_2), o(x_1, x_2))$.  
It induces 
a well-defined and bijective map 
\begin{equation}
\mu: \widetilde{Q}(F) \otimes \widetilde{Q}(F) \to \mathcal{C}_{\rm root}^{\rm s}(F).
\end{equation}
\end{thm} 

{\it Proof.}~ 
Consider  $(x_1 \cdot g, x_2 \cdot g)$ for $g \in F(\widetilde{Q}(F))$.  Then 
$x_1 \cdot g =[(D_1, a_1 \gamma)]$ and $x_2 \cdot g =[(D_2, a_2 \gamma)]$ for some loop $\gamma$ in $\R^4$ such that $\gamma$ misses $F$ and $\gamma(0)=\gamma(1)=q$.  
(This is seen from the definition of the operations $\ast$ and $\overline{\ast}$ for the knot symmetric quandle $\widetilde{Q}(F)$. See Figure~\ref{fgda.eps}.)  
The chord associated with 
$(x_1 \cdot g, x_2 \cdot g)$ is $\tilde{a_1} \gamma \gamma^{-1}\tilde{a_2}^{-1}$ 
which is homotopic to $\tilde{a_1}\tilde{a_2}^{-1}$ without changing $D_1$ and $D_2$.   Thus 
$(a(x_1 \cdot g, x_2 \cdot g),   o(x_1 \cdot g, x_2 \cdot g))$ is strongly equivalent to 
$(a(x_1, x_2), o(x_1, x_2))$.    Hence, $\mu$ is well-defined.  

Let $(a, o)$ be a root-oriented chord attached to $F$. Let $D_1$ (or $D_2$) be an oriented meridian disk of $F$ over the initial point $a(0)$ (or the terminal point $a(1)$) such that the orientation matches the root orientation $o$.  Moving $a$ up to homotopy, we change $a$ such that $a(1/2)=q$.  Then $a$ is the associated chord with $(x_1, x_2)$ for some 
$(x_1, x_2) \in \widetilde{Q}(F) \times \widetilde{Q}(F)$.  Thus, $\mu$ is surjective.  

Suppose that $[(a(x_1, x_2), o(x_1, x_2))] = [(a(x'_1, x'_2), o(x'_1, x'_2))]$.  
Let $\alpha_s: [0,1] \to \R^4$ $(s \in [0,1])$ be a homotopy of chords attached to $F$ 
such that $\alpha_0= a(x_1, x_2)$ and $\alpha_1= a(x'_1, x'_2)$.  
Let $\beta: [0,1] \to \R^4$ be the path defined by 
$\beta(s) = \alpha_s(1/2)$ for $s \in [0,1]$.  Note that $\beta$ misses $F$ and 
$\beta(0)=\beta(1)=q$.  Then $\tilde{a_1}\beta$ is homotopic to $\tilde{a'_1}$ 
through paths in $\R^4$ such that the paths intersect with $F$ at their initial points only and the teminal points are $q$.  This implies that $x_1\cdot g = x'_1 \in \widetilde{Q}(F)$ where $g$ is an element of $F(\widetilde{Q}(F))$ which presents  $[\beta] \in \pi_1(\R^4 \setminus F, q) = As(\widetilde{Q}(F), \rho)$.  (Note that there is a natural projection $F(\widetilde{Q}(F)) \to As(\widetilde{Q}(F), \rho)$ and $As(\widetilde{Q}(F), \rho)$ is identified with $G(F)= \pi_1(\R^4 \setminus F, q)$ in the sense of Theorem 4.5 of \cite{Kamada14A} (Remark~\ref{remark:group}), which is analogous to the fact given in \cite{FennRourke92} that $As(Q(F))$ is identified with $G(F)$.)  
Similarly, we have $x_2\cdot g = x'_2$ for the same $g$.  Thus $[(x_1, x_2)] = [(x'_1, x'_2)]$ in $\widetilde{Q}(F) \otimes \widetilde{Q}(F)$.  
\qed

\vspace{2mm} 

Combining 
the bijection 
$\phi :  \mathcal{H}^{\rm s}(F)  \to \mathcal{C}_{\rm root}^{\rm s}(F)$ in Theorem~\ref{thm:HCsw} and the inverse map of $\mu$ in Theorem~\ref{thm:QCs}, 
we obtain a bijection as follows. 

\begin{thm}\label{thm:HQs} 
Let $F$ be a surface-link. The map 
\begin{equation}
\Phi= \mu^{-1} \circ \phi:  \mathcal{H}^{\rm s}(F) \to \widetilde{Q}(F) \otimes \widetilde{Q}(F)
\end{equation}
is a bijection. 
\end{thm} 

{\it Proof.}~ 
It follows from Theorems~\ref{thm:HCsw} and \ref{thm:QCs}. \qed 

\vspace{2mm}

This theorem gives the former assertion of Theorem~\ref{thm:mainA}. 

We consider the case where $F$ is an oriented surface-link.

Let $x_1 = [(D_1, a_1)]$ and $x_2 = [(D_2, a_2)]$ be elements of $Q(F)$.  
Let $\tilde{a_1}$ and $\tilde{a_2}$ be extended paths of $a_1$ and $a_2$.  
The concatenation of $\tilde{a_1}$ and $\tilde{a_2}^{-1}$ is 
a chord attached to $F$.  We denote it by $a(x_1, x_2)$ and call it the {\it chord associated with } 
$(x_1, x_2)$.  

\begin{thm}\label{thm:QCsori} 
Let 
\begin{equation}
Q(F) \times Q(F) \to \mathcal{C}^{\rm s}(F).
\end{equation}
be the map sending $(x_1, x_2)$ to 
the strong equivalence class of $a(x_1, x_2)$.  It induces 
a well-defined and bijective map 
\begin{equation}
\mu_0^{\rm ori}: Q(F) \otimes Q(F) \to \mathcal{C}^{\rm s}(F).
\end{equation}
\end{thm} 

{\it Proof.}~ 
It is proved by a similar argument to the proof of Theorem~\ref{thm:QCs} by forgetting the root orientations. \qed 

\vspace{2mm} 

Combining the bijection 
$\phi_0^{\rm ori} :  \mathcal{H}^{\rm s}(F)^{\rm ori}  \to \mathcal{C}^{\rm s}(F)$ in Theorem~\ref{thm:A1} and the inverse map of $\mu_0^{\rm ori}$, we have a bijection as follows.  

\begin{thm}\label{thm:HQsori} 
Let $F$ be an oriented surface-link. The map 
\begin{equation}
\Phi_0^{\rm ori} = (\mu_0^{\rm ori})^{-1} \circ \phi_0^{\rm ori}:  \mathcal{H}^{\rm s}(F)^{\rm ori}  \to Q(F) \otimes Q(F)
\end{equation}
is a bijection. 
\end{thm} 

{\it Proof.}~ 
It follows from Theorems~\ref{thm:A1} and \ref{thm:QCsori}. \qed 

\vspace{2mm}

This theorem gives the former assertion of Theorem~\ref{thm:mainB}.

\begin{remark} 
Let $F$ be an oriented surface-link.  
There is a commutative diagram as follows: 
\def\spmapright#1{\smash{%
  \mathop{\hbox to 1.3cm{\rightarrowfill}} 
    \limits^{#1}}} 
\def\sbmapright#1{\smash{%
  \mathop{\hbox to 1.3cm{\rightarrowfill}} 
    \limits^{#1}}} 
\def\lmapdown#1{\Big\downarrow %
  \llap{$\vcenter{\hbox{$\scriptstyle#1\, $}}$ }} 
\def\rmapdown#1{\Big\downarrow %
  \rlap{$\vcenter{\hbox{$\scriptstyle#1\, $}}$ }} 
\begin{equation}
\begin{array}{ccccc} 
\mathcal{H}^{\rm s}(F)^{\rm ori} &\spmapright{\phi_0^{\rm ori}} & \mathcal{C}^{\rm s}(F) & \spmapright{(\mu_0^{\rm ori})^{-1}} & Q(F) \otimes Q(F) \\ 
\lmapdown{\rm incl.}  &   & \lmapdown{}   & &  \rmapdown{\rm incl.}  \\ 
\mathcal{H}^{\rm s}(F)              & \spmapright{\phi} & \mathcal{C}_{\rm root}^{\rm s}(F)  & \spmapright{\mu^{-1}} & \widetilde{Q}(F) \otimes \widetilde{Q}(F),
\end{array} 
\end{equation}
where the first and the third vertical maps are inclusion maps, and the second vertical map 
$\mathcal{C}^{\rm s}(F) \to \mathcal{C}_{\rm root}^{\rm s}(F)$ is an  injective map 
sending $[a]$ to $[(a, o_F)]$ where $o_F$ is the root orientation of $a$ which matches the orientation of $F$. 
All horizontal maps are bijections. 
\end{remark}


\section{$1$-handles and chords up to weak equivalence} \label{sect:weak}


In this section, we consider $1$-handles and chords  up to weak equivalence. 

Let $F$ be a surface-link and let $(\widetilde{Q}(F), \rho)$ be the knot symmetric quandle.  Let $\tau$ and $\rho$ be involutions of $\widetilde{Q}(F) \otimes \widetilde{Q}(F)$  with $\tau [x_1, x_2] = [x_2, x_1]$ and $\rho [x_1, x_2] = [\rho(x_1), \rho(x_2)]$.  

\begin{thm}\label{thm:QCw}
The bijection 
$\mu: \widetilde{Q}(F) \otimes \widetilde{Q}(F) \to \mathcal{C}_{\rm root}^{\rm s}(F)$ 
in Theorem~\ref{thm:QCs} induces a bijection 
\begin{equation}
\mu^{\langle \tau, \rho \rangle}: \widetilde{Q}(F) \otimes \widetilde{Q}(F)/ \langle \tau, \rho \rangle \to \mathcal{C}_{\rm root}^{\rm w}(F). 
\end{equation}
\end{thm}

{\it Proof.}~ 
Let ${\rm rev}_1$ and ${\rm rev}_2$ be involutions of 
$\mathcal{C}_{\rm root}^{\rm s}(F)$  
with ${\rm rev}_1([(a, o)]) = [({\rm rev}(a), o)]$ and  ${\rm rev}_2([(a, o)]) = [(a, -o)]$.  
Note that the quotient of $\mathcal{C}_{\rm root}^{\rm s}(F)$ 
by the equivalence relation generated by ${\rm rev}_1$ and ${\rm rev}_2$ is  $\mathcal{C}_{\rm root}^{\rm w}(F)$.  
Since $\mu \circ \tau = {\rm rev}_1 \circ \mu$ and $\mu \circ \rho = {\rm rev}_2 \circ \mu$, we have the result.  
\qed 

\vspace{2mm}

Combining  the bijection 
$\psi :  \mathcal{H}^{\rm w}(F)  \to \mathcal{C}_{\rm root}^{\rm w}(F)$ in Theorem~\ref{thm:HCsw} and the inverse map of $\mu^{\langle \tau, \rho \rangle}$,  we have a bijection as follows. 

\begin{thm}\label{thm:HQw} 
Let $F$ be a surface-link.  The map 
\begin{equation}
\Psi= (\mu^{\langle \tau, \rho \rangle})^{-1} \circ \psi:  \mathcal{H}^{\rm w}(F) \to \widetilde{Q}(F) \otimes \widetilde{Q}(F)/ \langle \tau, \rho \rangle
\end{equation}
is a bijection. 
\end{thm} 

{\it Proof.}~ 
It follows from Theorems~\ref{thm:HCsw} and \ref{thm:QCw}. 
\qed 

\vspace{2mm} 

This theorem give the latter assertion of Theorem~\ref{thm:mainA}.

Let $F$ be an oriented surface-link. 
Let $\tau$ be the involution of $Q(F) \otimes Q(F)$ 
with $\tau [x_1, x_2] = [x_2, x_1]$.   

\begin{thm}\label{thm:QCwori}
The bijection 
$\mu_0^{\rm ori}: Q(F) \otimes Q(F) \to \mathcal{C}^{\rm s}(F)$ induces a bijection 
\begin{equation}
\mu_0^{\langle \tau \rangle} : Q(F) \otimes Q(F)/{\langle \tau \rangle} \to \mathcal{C}^{\rm w}(F).  
\end{equation}
\end{thm} 

{\it Proof.}~ 
Let ${\rm rev}: \mathcal{C}^{\rm s}(F) \to \mathcal{C}^{\rm s}(F)$ be the involution of $Q(F)$  
with ${\rm rev}([a]) = [{\rm rev}(a)]$. 
The quotient  of 
$\mathcal{C}^{\rm s}(F)$ by the equivalence relation generated by ${\rm rev}$ is $\mathcal{C}^{\rm w}(F)$.  
Since $\mu_0 \circ \tau = {\rm rev} \circ \mu_0$, we have the result. \qed 

\vspace{2mm}

Combining the bijection 
$\psi_0^{\rm ori} :  \mathcal{H}^{\rm w}(F)^{\rm ori}  \to \mathcal{C}^{\rm w}(F)$ in Theorem~\ref{thm:A1} and the inverse map of $\mu_0^{\langle \tau \rangle}$, 
we have a bijection as follows. 

\begin{thm}\label{thm:HQwori} 
Let $F$ be an oriented surface-link.  The map 
\begin{equation}
\Psi_0^{\rm ori} = (\mu_0^{\langle \tau \rangle})^{-1} \circ \psi_0^{\rm ori}:  
\mathcal{H}^{\rm w}(F)^{\rm ori}  \to Q(F) \otimes Q(F)/ \langle \tau \rangle
\end{equation}
is a bijection. 
\end{thm} 

{\it Proof.}~ 
It follows from Theorems~\ref{thm:HCsw} and \ref{thm:QCwori}. \qed 

\vspace{2mm}

This theorem give the latter assertion of Theorem~\ref{thm:mainB}.


\section{Examples and invariants of $1$-handles} \label{sect:example}


\begin{example} 
Let $T(2,n)$ be a torus knot of type $(2,n)$ with $n=2m+1$, and 
let $F$ be the $2$-twist spun $T(2,n)$.   
The knot quandle $Q(F)$ has a presentation 
\begin{equation} 
\langle a, \, b \mid  b = a^{(ba)^m}, \, b^{a^2}=b \rangle, 
\end{equation} 
which 
 is isomorphic  to the dihedral quandle $R_n$ of order $n$ 
(cf. \cite{Satoh02}), where we follow the convention of quandle operation and presentation due to \cite{FennRourke92} (cf. \cite{KamadaBook2017}). 
The full knot quandle $\widetilde{Q}(F)$, with the standard involution as a knot symmetric quandle, is  isomorphic  to the symmetric double of $Q(F)$ (Remark~\ref{remark:double}) and hence to the symmetric double $D(R_n)$.   
Fix an isomorphism from $Q(F)$ to $R_n$ and identify $Q(F)$ with $R_n$.  Then 
\begin{equation} 
Q(F) \otimes Q(F) = R_n \otimes R_n   \quad \mbox{and} \quad 
\widetilde{Q}(F) \otimes \widetilde{Q}(F) = D(R_n) \otimes D(R_n).    
\end{equation} 

By Theorems~\ref{thm:HQsori} and~\ref{thm:HQs}, we obtain bijections 
$\mathcal{H}^{\rm s}(F)^{\rm ori} \to R_n \otimes R_n$ and 
$\mathcal{H}^{\rm s}(F) \to D(R_n) \otimes D(R_n)$.  
Since we have Theorems~\ref{thm:dihedral:odd} and \ref{thm:d:dihedral:odd}, 
they provide us complete classifications of oriented $1$-handles and all $1$-handles attached to $F$ up to strong equivalence.  In particular, 
we see that $\# \mathcal{H}^{\rm s}(F)^{\rm ori} = m+1$  and 
$\# \mathcal{H}^{\rm s}(F) = 2n+2$.   

By Theorems~\ref{thm:HQwori} and~\ref{thm:HQw}, we obtain bijections 
$\mathcal{H}^{\rm w}(F)^{\rm ori} \to R_n \otimes R_n/ \langle \tau \rangle$ and 
$\mathcal{H}^{\rm w}(F) \to D(R_n) \otimes D(R_n)/ \langle \tau, \rho \rangle$.  
Since we have Corollaries~\ref{cor:dihedral:odd} and \ref{cor:d:dihedral:odd:taurho}, 
they provide us complete classifications of oriented $1$-handles and all $1$-handles  up to weak equivalence.  In particular, 
we have $\# \mathcal{H}^{\rm w}(F)^{\rm ori} = m+1$  and 
$\# \mathcal{H}^{\rm w}(F) = n+1$.   
\end{example} 

\begin{example}
Let $F$ be a knotted projective plane which is the connected sum of  
the $2$-twist spun $T(2,n)$ with $n=2m+1$ and a (positive or negative) standard projective plane in $\R^4$.  The full knot quandle $\widetilde{Q}(F)$  is  isomorphic  to the  
dihedral quandle $R_n$ and  the standard involution $\rho$ as a knot symmetric quandle is the identity map (cf. \cite{Kamada14A, Satoh02}).   
Fix an identification of $\widetilde{Q}(F)$ and $R_n$, and we have 
\begin{equation} 
\widetilde{Q}(F) \otimes \widetilde{Q}(F) = R_n \times R_n.  
\end{equation}

By Theorem~\ref{thm:HQs}, we obtain a bijection 
$\mathcal{H}^{\rm s}(F) \to R_n \otimes R_n$, which provides us a complete classification of $1$-handles attached to $F$ up to strong equivalence.  In particular, we have 
$\# \mathcal{H}^{\rm s}(F) = m+1$. 

By Theorem~\ref{thm:HQw}, we obtain a bijection 
$\mathcal{H}^{\rm w}(F) \to R_n \otimes R_n / \langle \tau, \rho  \rangle = R_n \otimes R_n / \langle \tau \rangle$, which provides us a complete classification of $1$-handles attached to $F$ up to weak equivalence.  In particular, we have 
$\# \mathcal{H}^{\rm w}(F) = m+1$. 
\end{example}

For a surface-link such that  $\widetilde{Q}(F)$ (or $Q(F)$) is infinite or has a large cardinality, it is difficult to compute or list the elements of the tensor product $\widetilde{Q}(F) \otimes \widetilde{Q}(F)$ (or $Q(F) \times Q(F)$) and its quotient.  In such a case, it is useful to construct an invariant.

Let $F$ be a surface-link and let $X$ be a (finite) quandle.  Let 
$f : \widetilde{Q}(F) \to X$ be a quandle homomorphism.  
By Lemma~\ref{lem:tensorhom}, 
it induces a map 
\begin{equation}
f \otimes f : \widetilde{Q}(F) \otimes \widetilde{Q}(F) \to X \otimes X, \quad [x_1, x_2] \mapsto [f(x_1), f(x_2)].  
\end{equation}
Combining the bijection in Theorem~\ref{thm:HQs}, we have a map 
\begin{equation}
{\rm I}^{\rm s}(f) =  (f \otimes f) \circ \Phi 
 : \mathcal{H}^{\rm s}(F) \to X \otimes X,  \quad [h] \mapsto (f \otimes f)[\Phi([h])].  
\end{equation}
In this way, for each quandle homomorphism $f : \widetilde{Q}(F) \to X$, we obtain an invariant ${\rm I}^{\rm s}(f)$ of strong equivalence classes of $1$-handles attached to $F$ valued in $X \otimes X$. 

When $F$ is an oriented surface-link, for any quandle homomorphism $f : Q(F) \to X$, we obtain an invariant ${\rm I}^{\rm s}_{\rm ori}(f) = (f \otimes f) \circ \Phi_0^{\rm ori}$ of strong equivalence classes of oriented $1$-handles attached to $F$ valued in $X \otimes X$. 

Let $(X, \rho)$ be a symmetric quandle.  Let 
$f : (\widetilde{Q}(F), \rho) \to (X, \rho) $ be a symmetric quandle homomorphism to a (finite) symmetric quandle $(X, \rho)$.  
By Theorem~\ref{thm:HQw}, we have a map 
\begin{equation}
{\rm I}^{\rm w}(f)  = (f \otimes f) \circ \Psi 
: \mathcal{H}^{\rm w}(F) \to X \otimes X / \langle \tau, \rho \rangle ,  \quad [h] \mapsto (f \otimes f)[\Psi([h])].  
\end{equation}
In this way, for each symmetric quandle homomorphism $f : (\widetilde{Q}(F), \rho) \to (X, \rho) $, we obtain an invariant ${\rm I}^{\rm w}(f)$ of weak equivalence classes of $1$-handles attached to $F$ valued in $X \otimes X / \langle \tau, \rho \rangle$. 

When $F$ is an oriented surface-link, for any quandle homomorphism $f : Q(F) \to X$, we obtain an invariant ${\rm I}^{\rm w}_{\rm ori}(f) = (f \otimes f) \circ \Psi_0^{\rm ori}$ of weak equivalence classes of oriented $1$-handles attached to $F$ valued in $X \otimes X / \langle \tau \rangle$. 

Now we have Theorems~\ref{thm:mainC} and \ref{thm:mainD}.  

The dihedral quandle $R_n$ and the symmetric double $(D(X), \rho)$ can be used to construct such an invariant.   
Refer to \cite{Hayashi, LopesRoseman, Vendramin} for other finite quandles with small orders.  
For such a small quandle, one can compute the tensor product and use it to construct an invariant of $1$-handles.  

\vspace{0.4in}

\end{document}